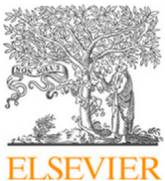
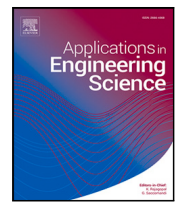

# Phantom Domain Finite Element Method: A novel approach for heterogeneous materials

Tianlong He, Philippe Karamian-Surville *, Daniel Choï

*Normandie Univ., UNICAEN, CNRS Laboratoire de Mathématiques Nicolas Oresme, 14000 Caen, France*

ABSTRACT

In this paper, we introduce the Phantom Domain Finite Element Method (PDFEM), a novel computational approach tailored for the efficient analysis of heterogeneous and composite materials. Inspired by fictitious domain methods, this method employs a structured mesh to discretize the entire material domain while utilizing separate, independent meshes for the inclusions. These inclusion meshes are coupled to the structured mesh via a substitution matrix, enabling them to act as phantom meshes that do not directly contribute to the final system of equations. This framework offers significant advantages, including enhanced flexibility in handling complex inclusion geometries and improved computational efficiency. To assess the accuracy and robustness of the proposed method, numerical experiments are conducted on structures containing inclusions of various geometries. In order to emphasize the efficiency of the PDFEM method, a numerical simulation is presented to highlight its advantages in the case of long natural fibers, such as flax and linen. These simulations are compared against FEM calculations, demonstrating the efficiency of PDFEM. Indeed, meshing such fine structures requires an extremely high number of elements, and in some cases, meshing becomes particularly challenging due to the complexity of the geometries.

## 1. Introduction

A composite material is a heterogeneous system consisting of at least two distinct constituents: a matrix, which typically acts as a binder, and inclusions composed of one or more materials. Inclusions are often referred to as reinforcements, as the primary objective of manufacturing composite materials is to integrate the properties of both the matrix and the inclusions, thereby yielding a material with enhanced characteristics that differ from those of the individual components. To efficiently model such materials, numerical homogenization techniques are employed to compute their effective properties.

To accurately model a domain made of heterogeneous material, we typically require a conforming mesh for Finite Element based Methods (Debordes, 1986) or a high-resolution image for iterative FFT-based methods (Moulinec and Suquet, 1994, 1995).

Finite Element Method (FEM) requires accurate calculation of integrals that define the energy bilinear form involved in the weak formulation of the model (Debordes, 1986; Bensoussan et al., 1978; Sanchez-Palencia and Sanchez-Hubert, 1992). This generally necessitates a conforming mesh. However, constructing such meshes for complex geometries is nontrivial (Schneider et al., 2016), and thin-shaped inclusions such as fibers can lead to very fine meshes, significantly increasing computational complexity.

To address these meshing challenges, several advanced approaches have been proposed, including the eXtended Finite Element Method (XFEM) (Legrain et al., 2016; Moës et al., 1999; Sukumar et al., 2001) and the Finite Cell Method (FCM) (Düster et al., 2008; Parvizian et al., 2007; Duczek and Gabbert, 2016; Kudela et al., 2020).

The XFEM integrates the Level Set Method to define inclusions while incorporating enriched shape functions, thereby enhancing its capability to model complex geometries (Yvonnet et al., 2008). The FCM employs a single structured mesh, where the elements are referred to as cells, and improves the accuracy of integral approximations through cell subdivisions (Parvizian et al., 2007; Düster et al., 2008, 2012; Legrain et al., 2016).

The core concept of the Phantom Domain Finite Element Method (PDFEM), introduced in this paper, is similarly to approximate the integrals that define the energy bilinear form using independent meshes for the inclusions, while establishing a relationship with the *structured mesh* of the overall domain. PDFEM is inspired by Fictitious Domain Methods (FDM); however, unlike traditional FDM, it does not employ Lagrange multipliers to enforce coupling between the meshes (Glowinski et al., 1994; Haslinger and Renard, 2009).

In PDFEM, integrals are computed over both the entire domain and the inclusions using independent meshes. The degrees of freedom of the inclusion meshes are then mapped to those of the structured domain mesh via a substitution matrix, which can be efficiently computed. The PDFEM, as formulated in this work, is directly applicable to inclusions of arbitrary shape or geometry, provided that their corresponding meshes are available.

The paper is structured as follows: after this introduction, we present the Phantom Domain Finite Element Method (PDFEM) for a composite material consisting of a matrix and inclusions. For clarity, the method is first described for linear static thermal problems, formulated as boundary value problems. The presentation of the PDFEM is organized into three key steps:

- Splitting the energy functional
- Constructing the substitution matrix
- Ensuring sufficient pixelization of inclusions

Numerical experiments are then conducted for both thermal and elastic boundary value problems. The results obtained using PDFEM are compared against those from conforming Finite Element Methods (FEM) for various inclusion geometries.

Finally, we extend our study to include a random fiber distribution to assess the method's capability in computing global material properties.

*Notations*

In this paper, scalar quantities are represented in *italic*, while vector quantities are denoted in **boldface**. For instance, the *temperature* is expressed as the scalar $u$, whereas the *displacement* is described as the vector field $\mathbf{u}$. The symbols div and $\nabla$ correspond to the divergence and gradient differential operators, respectively. Furthermore, $\hat{u}$ represents an element of the Euclidean space $\mathbb{R}^n$.

## 2. The phantom domain finite element method (PDFEM)

For the sake of simplicity, the PDFEM for thermal problems is presented in detail. Let $\Omega$ be a domain constituted by a heterogeneous material as shown in Fig. 1. Under the action of an external heat source $f$, imposed temperature $u_0$ on $\Gamma_D$ and imposed heat flux $F$ on $\Gamma_N$, the thermal equilibrium boundary value problem is written as:

$$\begin{cases} \operatorname{div} \mathbf{q} + f = 0 & \text{in } \Omega \\ u = u_0 & \text{on } \Gamma_D \\ \mathbf{q} \cdot \mathbf{n} = F & \text{on } \Gamma_N, \end{cases} \quad (1)$$

where $\Gamma_D \cup \Gamma_N = \partial \Omega$ and $\Gamma_D \cap \Gamma_N = \emptyset$, $u$ denotes the temperature, and $\mathbf{q}$ is the heat flux.

Theoretically, $\Omega$ can have any arbitrary shape, but for practical implementation, we consider domains such as $\Omega = [0,1]^d$ with $d = 2, 3$ since we wish to use a structured mesh. The thermal conductivity tensor $\Lambda$ is defined as:

$$\mathbf{q} = \Lambda \nabla u. \quad (2)$$

For simplicity, let the material domain $\Omega$ consist of two homogeneous and isotropic media: the matrix $\Omega_{\text{mat}} = \Omega \setminus \Omega_{\text{inc}}$ and inclusions $\Omega_{\text{inc}}$. We assume the conductivity matrix $\Lambda$ to be piecewise constant:

$$\Lambda(x) = \begin{cases} \Lambda_{\text{mat}}, & x \in \Omega_{\text{mat}}, \\ \Lambda_{\text{inc}}, & x \in \Omega_{\text{inc}}. \end{cases}$$

The variational problem equivalent to the boundary value problem (1) is:

$$\begin{cases} \text{Find } u \in V_{\text{adm}} = \{u \in H^1(\Omega) \mid u|_{\Gamma_D} = u_0\} \\ a(u, u^*) = L(u^*) \quad \forall u^* \in V_0 \end{cases} \quad (3)$$

where $V_0$ is the vector subspace associated with the (affine) admissible space $V_{\text{adm}}$:

$$V_0 = \{u \in H^1(\Omega) \mid u|_{\Gamma_D} = 0\}.$$

and where

$$a(u, u^*) = \int_\Omega (\Lambda \nabla u) \cdot \nabla u^*, \qquad L(u^*) = \int_\Omega f u^* + \int_{\Gamma_N} F u^*.$$

Classically, the problem (3) is well-posed within the Lax–Milgram framework (Rektorys, 1977). The unique solution to the problem (3) also minimizes the energy functional $J(u)$ in $V_{\text{adm}}$:

$$J(u) = \frac{1}{2} a(u, u) - L(u). \quad (4)$$

We define an interpolation space $V^h$ of dimension $n$ such that any $u^h \in V^h$ can be represented by $\hat{u} \in \mathbb{R}^n$. We define the stiffness (conductivity) matrix $K$ and vector $L$ as:

$$J(u^h) = \frac{1}{2} a(u^h, u^h) - L(u^h), \quad (5)$$

$$J(u^h) = \frac{1}{2} \hat{u}^\top K \hat{u} - \hat{u}^\top L. \quad (6)$$

The basic idea of the Phantom Domain Finite Element Method (PDFEM) is based on three main steps:

a – *Splitting the Energy Form $J$*: The energy functional $J$ is split into two components, one representing the entire domain and another for the inclusions.

b – *Calculation of Matrices Defined by the Split*: These matrices are calculated using independent meshes for the domain and the inclusions.

c – *Construction of a Substitution Matrix*: This matrix links the degrees of freedom (DOFs) of the inclusion mesh to those of the entire mesh.

*Splitting of the energy form $J$*

Since the material is composed of two constituents, the energy functional (4) is split into corresponding parts:

$$J(u) = \frac{1}{2} \int_{\Omega \setminus \Omega_{\text{inc}}} (\Lambda_{\text{mat}} \nabla u) \cdot \nabla u + \frac{1}{2} \int_{\Omega_{\text{inc}}} (\Lambda_{\text{inc}} \nabla u) \cdot \nabla u \quad - L(u). \quad (7)$$

The main idea of PDFEM is to use distinct and *a priori* incompatible and independent meshes representing the whole domain $\Omega$ and the inclusions $\Omega_{\text{inc}}$ instead of a unique conformal mesh that matches the geometry of the inclusions. See Fig. 1 for an illustration in the case of one disk inclusion in $\Omega = [0, 1]^2$.

As we have:

$$\int_{\Omega \setminus \Omega_{\text{inc}}} (\Lambda_{\text{mat}} \nabla u) \cdot \nabla u = \int_\Omega (\Lambda_{\text{mat}} \nabla u) \cdot \nabla u - \int_{\Omega_{\text{inc}}} (\Lambda_{\text{mat}} \nabla u) \cdot \nabla u,$$

from (7), the energy functional can be rewritten as:

$$J(u) = \underbrace{\frac{1}{2} \int_\Omega (\Lambda_{\text{mat}} \nabla u) \cdot \nabla u}_{J_{\text{mat}}} + \underbrace{\frac{1}{2} \int_{\Omega_{\text{inc}}} ([\Lambda_{\text{inc}} - \Lambda_{\text{mat}}] \nabla u) \cdot \nabla u}_{J_{\text{inc}}} \quad - L(u). \quad (8)$$

The two split integrals $J_{\text{mat}}$ and $J_{\text{inc}}$ are defined over the whole domain $\Omega$ and the inclusion $\Omega_{\text{inc}}$, respectively.

The use of distinct meshes allows independent numerical computation of $J_{\text{mat}}$ and $J_{\text{inc}}$. The functional $J_{\text{mat}}$ can be calculated in $V^h$ with the help of an $n \times n$ matrix $K_{\text{mat}}$:

$$J_{\text{mat}}(u^h) = \frac{1}{2} \int_\Omega (\Lambda_{\text{mat}} \nabla u^h) \cdot \nabla u^h = \frac{1}{2} \hat{u}^\top K_{\text{mat}} \hat{u}. \quad (9)$$

With a conformal mesh of the inclusion $\Omega_{\text{inc}}$ independent of the mesh defining $\Omega$, we define the interpolation space $W^l$ of finite dimension $p$, associated with the inclusion's mesh. Thus, any $v^l \in W^l$ can be represented by a vector $\hat{v} \in \mathbb{R}^p$.

To simplify the presentation, the meshes are denoted as $\Omega$ and $\Omega_{\text{inc}}$. The functional $J_{\text{inc}}$ can be defined in $W^l$ with a $p \times p$ matrix $K_{\text{inc}}$:

$$J_{\text{inc}}(v^l) = \frac{1}{2} \int_{\Omega_{\text{inc}}} ([\Lambda_{\text{inc}} - \Lambda_{\text{mat}}] \nabla v^l) \cdot \nabla v^l = \frac{1}{2} \hat{v}^\top K_{\text{inc}} \hat{v}. \quad (10)$$

Both matrices $K_{\text{mat}}$ and $K_{\text{inc}}$ can be computed using standard finite element procedures.

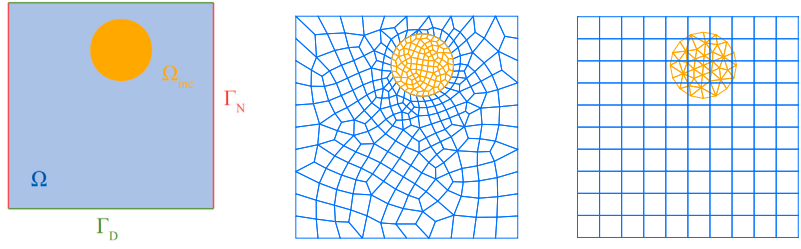

**Fig. 1.** A domain $\Omega$ with one inclusion $\Omega_{\text{inc}}$ (right) under the imposed temperature $u_0$ on $\Gamma_D$ and imposed heat flux $F$ on $\Gamma_N$, a conformal mesh (middle) for FEM and independent inclusion mesh and structured mesh (right) for PDFEM.

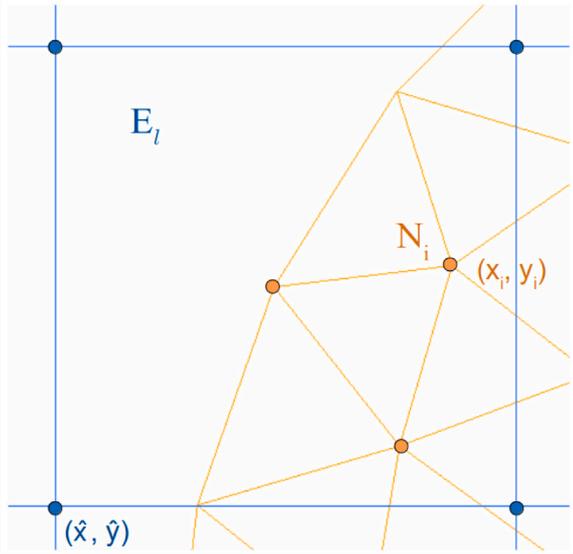

**Fig. 2.** A node $N_i$ from a triangular mesh of an inclusion in the quadrangular element $E_l$ of the mesh $\Omega$. Note that the nature of the inclusion's mesh (tri3) can be different from the structured mesh's (qua4).

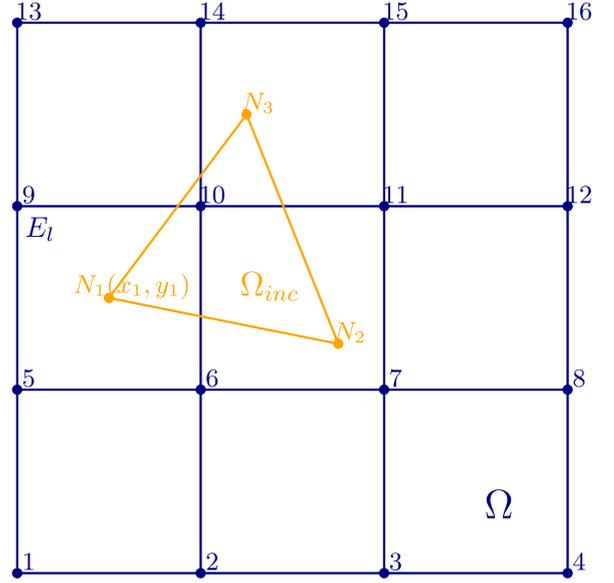

**Fig. 3.** A triangular mesh of an inclusion in the structured mesh $\Omega$ with mesh resolution equal to 3.

*Construction of a substitution matrix*

We defined the matrices $K_{\text{mat}}$ and $K_{\text{inc}}$, which allow us to separately compute the energy functionals $J_{\text{mat}}$ and $J_{\text{inc}}$. The objective here is to establish a relationship between $\hat{v}$, the degrees of freedom (DOFs) of the inclusion $\Omega_{\text{inc}}$, and $\hat{u}$, the DOFs of the structured mesh $\Omega$.

With a $n_x \times n_y$ structured mesh, we can define the elements of the mesh with the help of integer coordinates $(i, j)$ with $i = 1, \ldots, n_x$ and $j = 1, \ldots, n_y$. Given a node of coordinate $(x, y)$, it is straightforward to determine which element contains this node. Such a determination becomes nontrivial in the general case of non-structured meshes.

Thus, for each component $v_i$ of $\hat{v}$, we have $v_i = v(N_i)$, where $N_i$ is a node of $\Omega_{\text{inc}}$, let $E_l$ be an element of $\Omega$ in which node $N_i$ is included, as illustrated in Fig. 2 for a two-dimensional case.

Let $\hat{x}$, $\hat{y}$, and $\hat{z}$ be the coordinates of the nodes constituting the element $E_l$. In two dimensions, $\hat{z}$ is naturally set to 0. With an isoparametric element, we define a reference element and its corresponding shape functions $\hat{\Phi}$. Let $(x_i, y_i, z_i)$ be the coordinates of the node $N_i$, and let $(r_i, s_i, t_i)$ represent its coordinates in the reference element associated with $E_l$. By definition, we have:

$$x_i = \hat{\Phi}(r_i, s_i, t_i)^\top \hat{x},$$
$$y_i = \hat{\Phi}(r_i, s_i, t_i)^\top \hat{y},$$
$$z_i = \hat{\Phi}(r_i, s_i, t_i)^\top \hat{z}.$$

Let $v_i$ be the value of $\hat{v}$ at node $N_i$:

$$v_i = v(x_i, y_i, z_i). \tag{11}$$

Since we are working with isoparametric elements, we also have:

$$v(x_i, y_i, z_i) = \hat{\Phi}(r_i, s_i, t_i)^\top \hat{u}. \tag{12}$$

In other words, the value $v_i$ at a node of the inclusion's mesh is related to the values $\hat{u}$ on the mesh of the element $E_l$:

$$v_i = \hat{\Phi}(r_i, s_i, t_i)^\top \hat{u}. \tag{13}$$

This relation is extended (with zeros filled in) to all the values of $\hat{u}$ at the other nodes of the mesh $\Omega$. Let us denote $S_i$ as the $1 \times n$ line matrix, where $n$ is the total number of DOFs of the structured mesh $\Omega$, so that we have:

$$v_i = v(x_i, y_i, z_i) = \hat{\Phi}(r_i, s_i, t_i)^\top \hat{u} = S_i \hat{u},$$

where $\hat{u}$ (abusively) represents the nodal values of $u$ in an element $E_l$ but also in the entire domain $\Omega$.

To illustrate this, we provide an example in the two-dimensional case (which can be similarly extended to three dimensions). Consider a structured mesh of $\Omega$ with resolution 3 (i.e., $3 \times 3$ elements and $4 \times 4$ nodes), and a mesh of triangular inclusion $\Omega_{\text{inc}}$ with 1 element and 3 nodes, as shown in Fig. 3.

It can be seen that the node $N_1$ of $\Omega_{\text{inc}}$ is included in an element $E_l$ of $\Omega$ with nodes 5, 6, 9, 10.

We denote $(r_1, s_1)$ as the coordinates of node $N_1$ in the reference element associated with $E_l$. Within the isoparametric element $E_l$, we have the value $v_1$ at node $N_1$:

$$v_1 = v(x_1, y_1) = \hat{\Phi}(r_1, s_1)^\top \hat{u}.$$

The first line of the substitution matrix $S_1$, corresponding to node $N_1$ of the inclusion mesh $\Omega_{\text{inc}}$, is defined as:

$$[v_1] = \underbrace{\begin{bmatrix} 0 & \cdots & \Phi_1 & \Phi_2 & 0 & 0 & \Phi_4 & \Phi_3 & \cdots & 0 \end{bmatrix}}_{S_1} \begin{bmatrix} \vdots \\ u_5 \\ u_6 \\ u_7 \\ u_8 \\ u_9 \\ u_{10} \\ \vdots \end{bmatrix}.$$

The shape functions $\Phi_1$, $\Phi_2$, $\Phi_3$, and $\Phi_4$ correspond to the four nodes $u_5$, $u_6$, $u_{10}$, and $u_9$ of the element $E_l$ in the structured mesh $\Omega$.

The expression for $v_1$ indicates that it is a linear combination of shape functions. Specifically, it represents the value $v_1$ at node $N_1$ in terms of the shape functions $\hat{\Phi}(r_1, s_1)$ evaluated at the reference coordinates $(r_1, s_1)$, the coefficients being the degrees of freedom $\hat{u}$. Thus, $v_1$ is computed as:

$$v_1 = v(x_1, y_1) = \hat{\Phi}(r_1, s_1)^\top \hat{u}$$

This formulation shows that $v_1$ is the weighted sum of the nodal values, where the shape functions serve as the weights. This is a common approach in finite element analysis, where the field variable at a point within an element is determined by a linear combination of the nodal values, with the shape functions determining how these nodal values influence the field at that point.

Considering all the nodes of $\Omega_{\text{inc}}$, we finally obtain the substitution matrix $S$, constructed line by line, which substitutes $\hat{v}$ (the DOFs associated with the inclusion mesh $\Omega_{\text{inc}}$) for $\hat{u}$ (the DOFs associated with the mesh $\Omega$):

$$\hat{v} = S\hat{u}. \tag{14}$$

It is important to emphasize that the construction of the substitution matrix $S$ is fully formulated for structured meshes and does not require any sorting or testing steps during the procedure. Consequently, the construction of the substitution matrix $S$ is very cost effective.

We are now able to rewrite the quadratic energy form computed on the inclusion $J_{\text{inc}}$:

$$\begin{aligned} J_{\text{inc}} &= \frac{1}{2} \hat{v}^\top K_{\text{inc}} \hat{v} \\ &= \frac{1}{2} (S\hat{u})^\top K_{\text{inc}} S\hat{u} \\ &= \frac{1}{2} \hat{u}^\top \left( S^\top K_{\text{inc}} S \right) \hat{u}. \end{aligned}$$

Finally, we obtain the quadratic form of the main energy defined only on $\hat{u}$:

$$\begin{aligned} J(u) &= J_{\text{mat}}(u) + J_{\text{inc}}(u) - \hat{u}^\top L \\ &= \frac{1}{2} \hat{u}^\top K_{\text{mat}} \hat{u} + \frac{1}{2} \hat{u}^\top \left( S^\top K_{\text{inc}} S \right) \hat{u} - \hat{u}^\top L \\ &= \frac{1}{2} \hat{u}^\top \left( K_{\text{mat}} + S^\top K_{\text{inc}} S \right) \hat{u} - \hat{u}^\top L. \end{aligned} \tag{15}$$

In other words, with the matrix defined in (15), we apply a finite element method with the interpolation space defined on the structured mesh of domain $\Omega$.

**Remark.** The substitution matrix principle is also used to compare the PDFEM with the FEM as it enables the projection of the PDFEM solution into the FEM mesh.

### 2.1. Weak continuities at the interfaces

In the case of heterogeneous materials, weak continuity in the solution is expected at the material interface. As in PDFEM, the finite element calculation is performed only on the structured mesh of the domain. Consequently, it is evident that weak continuity cannot be maintained using standard isoparametric shape functions, leading to potential discrepancies in the solution at the matrix/inclusion interfaces.

However, since the method is based on energy minimization, these local discrepancies should not significantly impact the global response or overall properties of the system. This is because the energy minimization approach aims to determine the lowest possible total energy state for the entire structure rather than focusing solely on specific areas or interfaces. Therefore, while differences may arise at the matrix/inclusion interfaces, these discrepancies should not substantially affect the system's overall behavior.

While our approach is simple, efficient, and computationally fast, it is not suited for obtaining local results in the transition zone between materials. We acknowledge that our method is not designed for high-precision local field resolution and should not be used in applications requiring detailed stress distribution analysis at material interfaces. Instead, it serves as a computationally efficient tool for estimating effective material properties, offering a simpler alternative to more complex fictitious domain methods.

### 2.2. Higher order interpolation

Currently, PdFEM remains constrained to linear interpolation within the structured mesh. The use of higher-order standard polynomials can introduce oscillations in the substitution matrix, potentially leading to method instability (Joulaian and Düster, 2013). However, since the primary objective is to evaluate material properties, higher-order interpolation appears unnecessary and does not present a significant issue for global effective properties, as demonstrated by the observed linear convergence and numerical simulations conducted for various complex geometries. Therefore, this limitation does not represent a major obstacle for our method. Nevertheless, our approach should be further improved to accommodate higher-order interpolation. To achieve this, alternative shape functions must be explored.

### 2.3. Structured mesh and inclusions meshes for the PDFEM

Let us consider a $d$-dimensional domain $\Omega \in [0,1]^d$, where $d = 2, 3$, representing the medium of a composite material that may embed one or more inclusions. To implement the Phantom Finite Element Method, we require one structured mesh for the domain $\Omega$ and an independent mesh for each inclusion $\Omega_{\text{inc}}$. we require:

- One structured mesh for the domain $\Omega$
- An independent mesh for each inclusion $\Omega_{\text{inc}}$

The structured mesh of the domain $\Omega$ is defined by the coordinates of the nodes and their numbering through simple formulas. It is characterized by its resolution $n$. In the case of the square, a resolution $n$ corresponds to a $n \times n$ grid of elements, with $(n+1) \times (n+1)$ nodes when $d = 2$. The characteristic length of the elements is then given by $h = \frac{1}{n}$.

*Enveloping pixelization of the inclusions*

We define here the *pixelization* of an inclusion mesh, where each node of the inclusion mesh defines an element (containing the node) of the structured mesh. Such an element can be viewed as a *pixel* where the inclusion material is present. The union of all these pixels is referred to as the pixelization of the inclusion, as derived from its mesh.

An enveloping pixelization depends not only on the mesh of the inclusion but also on the resolution of the structured mesh. If the inclusion mesh is not fine enough, the pixelization might fail to be enveloping and could leave gaps. This is illustrated in Fig. 4, where we consider a single triangular inclusion in the domain $\Omega = [0,1]^2$. The inclusion mesh remains the same, but two different resolutions of the structured mesh are used. The pixels defined from the inclusion's mesh are shown in pink.

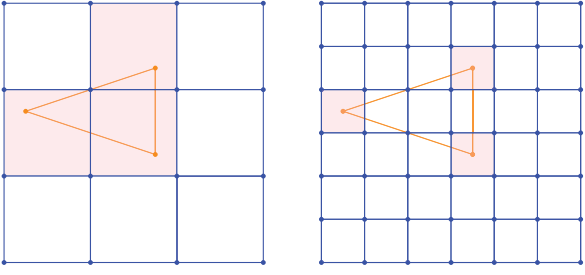

**Fig. 4.** Pixelization of a one-element mesh of a triangular inclusion with two different resolutions: $3 \times 3$ (left) and $6 \times 6$ (right) of the structured meshes of $\Omega = [0,1]^2$.

- For the $3 \times 3$ mesh, the pink elements fully envelop the inclusion
- For the $6 \times 6$ mesh, some gaps are visible.

In other words, the finer the structured mesh $\Omega$, the finer the inclusion's mesh $\Omega_{\text{inc}}$ must be to define an enveloping pixelization. Ensuring this condition is crucial for achieving an accurate approximation of the integrals. In the sequel, we will ensure that this condition is met :

In other words, as the structured mesh $\Omega$ becomes finer, the inclusion mesh $\Omega_{\text{inc}}$ must also be refined accordingly to achieve an accurate pixelization. Satisfying this condition is essential for obtaining precise integral approximations. In the following, we will ensure that this requirement is met.

To that purpose, we define a ratio $\eta$ as:

$$\eta = \frac{h_{\text{mat}}}{h_{\text{inc}}} \tag{16}$$

where $h_{\text{mat}}$ and $h_{\text{inc}}$ are the characteristic lengths of the meshes of the matrix $\Omega_{\text{mat}}$ and the inclusion $\Omega_{\text{inc}}$, respectively. A value of $\eta$ greater than 1 indicates that the inclusion mesh is finer than the matrix mesh, thereby ensuring an enveloping pixelization of the inclusion. Numerical experiments showed that a finer inclusion mesh does not significantly improve the effectiveness of the method.

## 3. Numerical experiments

We begin by presenting the specifications used for the numerical results of static thermal and elastic problems.

*Contrast parameter of the constitutive law*

For simplicity, we consider isotropic constitutive laws. We define a contrast parameter as the ratio between the characteristic coefficients of the inclusions and the matrix:

$$c_{\text{thermal}} = \frac{\lambda_{\text{inc}}}{\lambda_{\text{mat}}},$$

$$c_{\text{elastic}} = \frac{k_{\text{inc}}}{k_{\text{mat}}},$$

where $\lambda_{\text{inc}}$ (resp. $\lambda_{\text{mat}}$) denotes the thermal conductivity of the inclusion (resp. the matrix), and $k_{\text{inc}}$ (resp. $k_{\text{mat}}$) represents the bulk modulus of the inclusion (resp. the matrix). In our elastic tests, we assume an equal Poisson's ratio $\nu$ for both the inclusion and the matrix.

*SUBC and KUBC boundary conditions*

To compute effective properties, we consider pure Dirichlet and pure Neumann boundary conditions, referred to as the Kinematic Uniform Boundary Conditions (KUBC) and Static Uniform Boundary Conditions (SUBC) methods, respectively, for numerical homogenization (Kanit et al., 2003).

For simplicity, we present the boundary conditions in the two-dimensional case, with the extension to the three-dimensional case being straightforward.

- **Pure Dirichlet Boundary Value Problem (KUBC):**

    - **Thermal Case:** We impose a uniform temperature gradient at the boundary of $\Omega$, with $\mathbf{g} \in \mathbb{R}^2$ being a constant vector:

      $u(\mathbf{x}) = \mathbf{g} \cdot \mathbf{x}, \quad \forall \mathbf{x} \in \partial \Omega$

    - **Elastic Case:** We impose

      $\mathbf{u} = \overline{E}\,\mathbf{x}, \quad \forall \mathbf{x} \in \partial \Omega$

      where $\overline{E}$ is a constant $2 \times 2$ tensor, e.g., $\overline{E} = \begin{bmatrix} 1 & 0 \\ 0 & 0 \end{bmatrix}$ in the numerical example presented.

- **Pure Neumann Boundary Value Problem (SUBC):**

    - **Thermal Case:** We impose a uniform heat flux at the boundary of $\Omega$, with $\mathbf{Q} \in \mathbb{R}^2$ being a constant vector:

      $\mathbf{q}(\mathbf{x}) \cdot \mathbf{n} = \mathbf{Q} \cdot \mathbf{n}, \quad \forall \mathbf{x} \in \partial \Omega$

      where $\mathbf{n}$ denotes the outward-pointing unit normal at each point on the boundary $\partial \Omega$.

    - **Elastic Case:** We impose a uniform traction vector at the boundary:

      $\overline{\sigma}(\mathbf{x})\,\mathbf{n} = \overline{\Sigma}\,\mathbf{n}, \quad \forall \mathbf{x} \in \partial \Omega$

      where $\overline{\Sigma}$ is a constant $2 \times 2$ tensor, e.g., $\overline{\Sigma} = \begin{bmatrix} 0 & 0 \\ 0 & 1 \end{bmatrix}$ in the numerical example presented.

*Error estimation*

An a priori error estimate has not been defined. However, in some cases, we can estimate the difference between the PDFEM and standard FEM calculations. The relative difference between PDFEM and FEM can be considered as an a posteriori error estimate.

The comparison is made possible by the "projection" of the PDFEM solution $u^{\text{PDFEM}}$ onto an FEM conformal mesh. Specifically, we construct a substitution matrix $S$ in a manner similar to that in PDFEM, and the degrees of freedom $\hat{u}$ of the conformal mesh can be related to the degrees of freedom $\hat{u}^{\text{PDFEM}}$ defined on the structured mesh as follows:

$\hat{u} = S\hat{u}^{\text{PDFEM}}.$

This substitution matrix allows us to extrapolate the values $\hat{u}$ defined on the FEM conformal mesh from the values $\hat{u}^{\text{PDFEM}}$ calculated using PDFEM. We can then calculate the PDFEM/FEM relative difference of the solution between the two methods:

$$\text{relative difference} = \frac{\|\hat{u}^{\text{fem}} - S\hat{u}^{\text{PDFEM}}\|}{\|\hat{u}^{\text{fem}}\|}. \tag{17}$$

We also define the PDFEM/FEM relative difference using the $L^2$ norm and $H^1$ semi-norm:

$$\|u\|^2_{L^2(\Omega)} = \int_\Omega u^2,$$

$$|u|^2_{H^1(\Omega)} = \int_\Omega \|\nabla u\|^2.$$

In practice, the $L^2$ norm and $H^1$ semi-norm are computed with the help of the mass matrix $M$ and the stiffness matrix $K$.

*Computing resources*

The computations are performed using a finite element library code developed by the authors in Python/Fortran.

The meshes for the inclusions are generated with the Gmsh software (Geuzaine and Remacle, 2009).

All of the computations presented in this paper can be performed on personal computers (Intel Core i5, Core i7 and M1 ARM processor and 16 GB of memory were used for the simulations).

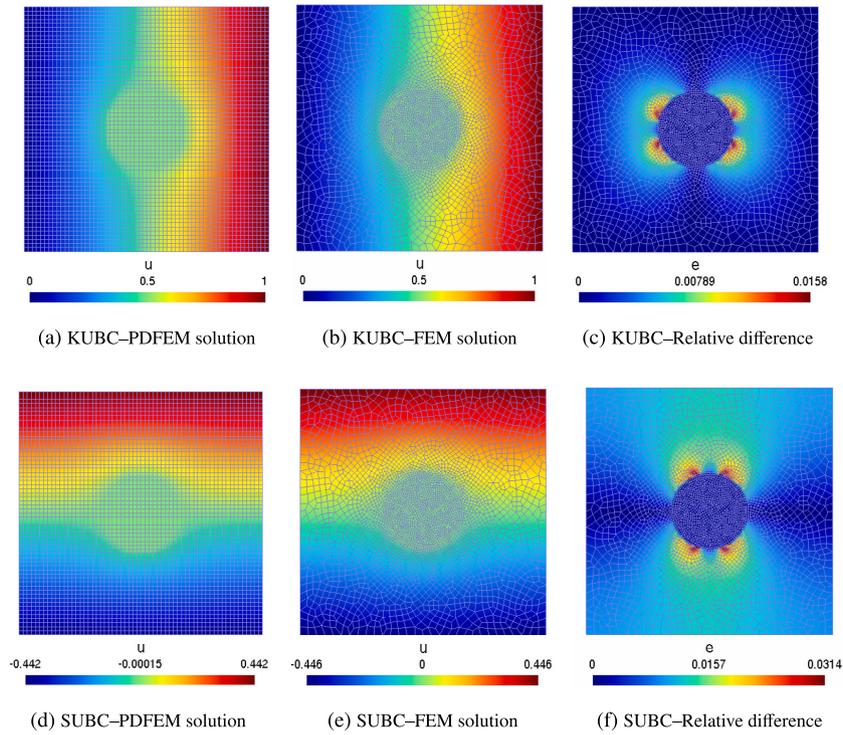

**Fig. 5.** Temperature computed by PDFEM, FEM, and the relative difference PDFEM/FEM for thermal problems: KUBC in (*a*), (*b*), (*c*) and SUBC in (*d*), (*e*), (*f*). Case of a disk inclusion with diameter $d = 0.3$.

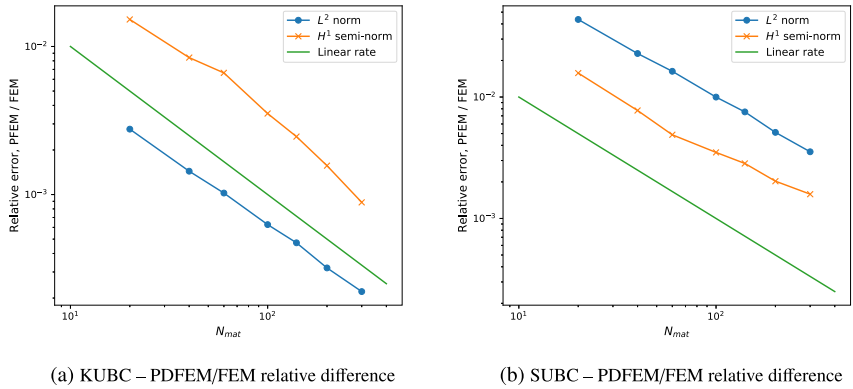

**Fig. 6.** PDFEM/FEM relative difference in $L^2$ norm and $H^1$ semi-norm for static thermal problems with respect to the mesh resolution. Case of a disk inclusion with $d = 0.3$.

*3.1. Case of a single inclusion*

We begin with the case of a single inclusion. We chose a disk or a square with respective diameter and side length of 0.3 in a $\Omega = [0, 1]^2$ square domain. The structured mesh of $\Omega$ has a resolution of $N_{\text{mat}} = 60$ (3721 nodes and 3600 quadrilateral elements). An independent mesh for the inclusion is defined for the PDFEM calculation. Additionally, a conformal mesh of the matrix-plus-inclusion domain is created to perform a classical FEM calculation, with the same order of number of nodes and elements. We then study the PDFEM/FEM relative difference, as defined in Eq. (17).

*Thermal problems*

In Figs. 5 and 7, temperatures are shown from PDFEM and FEM calculations, along with the PDFEM/FEM relative difference (as defined in Eq. (17)), for disk and square inclusions, respectively. In both the KUBC and SUBC boundary value problems, we observe that the maximum relative difference is of the order of $10^{-2}$ and is localized at the matrix/inclusion interface; it is more pronounced at the corners, in the case of a square inclusion. This is expected and is consistent with the fact that weak continuity (kink) is not satisfied at all in the matrix/inclusion interface by the interpolation polynomials from the structured mesh, see Section 2.1. A potential solution to overcome this issue might involve the introduction of modified shape functions within elements surrounding the interfaces.

To illustrate the evolution of the PDFEM/FEM relative difference with respect to the mesh resolution, we show in Fig. 6 the linear convergence of the PDFEM/FEM relative difference in both $L^2$ norm and $H^1$ semi-norm as the mesh resolution increases.

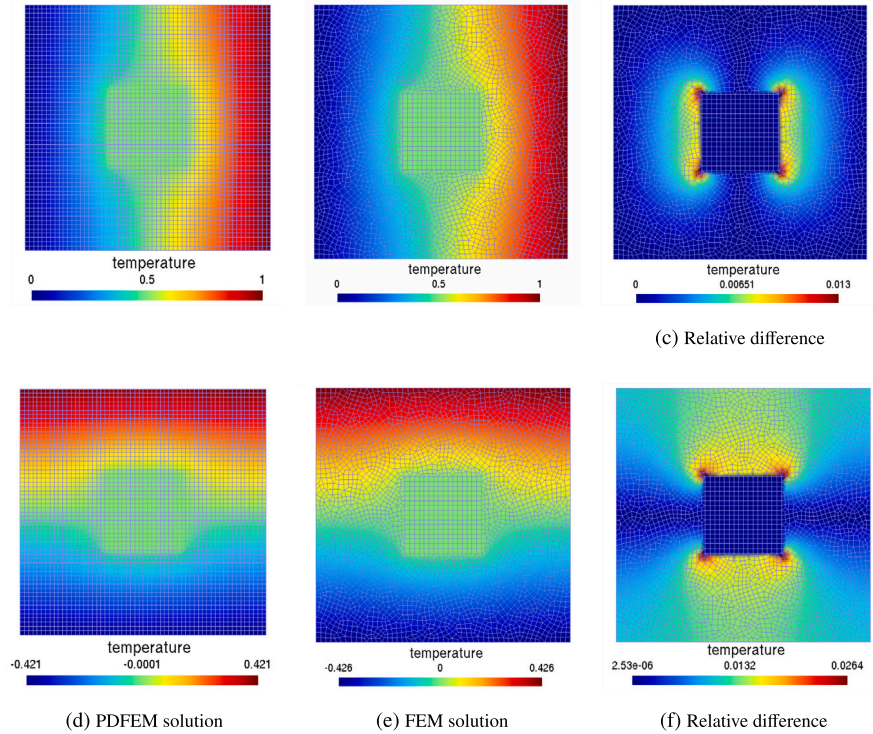

**Fig. 7.** Temperature computed by PDFEM, FEM, and the relative difference PDFEM/FEM for thermal problems: KUBC in (*a*), (*b*), (*c*) and SUBC in (*d*), (*e*), (*f*). Case of a square inclusion with side length $l = 0.3$.

*Elastic problems*

In the elastic case, we chose to display the PDFEM/FEM relative difference in displacements, strains, and Von-Mises equivalent stresses, as shown in Figs. 8 and 10.

The maximum relative difference in displacement is on the order of magnitude $10^{-2}$, while the maximum relative differences in strain and stress can approach approximately $\simeq 0.5$, as demonstrated in Fig. 8.

For a square inclusion, it is observed that the peak values of the relative difference are still concentrated at the four corners, with a maximum of $\simeq 0.02$ for displacement and $\simeq 0.2$ for both strain and Von-Mises stress. These differences are present as anticipated since they are localized at the matrix/inclusion interfaces due to the lack of weak continuity at the interface, as discussed in Section 2.1.

As shown in Fig. 9, linear convergence with respect to the resolution is observed for both the $L^2$ norm and the $H^1$ semi-norm of the PDFEM/FEM relative difference, with only a slight deviation in the latter.

*Special case of matching meshes*

In certain cases, such as with a square inclusion, the nodes and edges of the square's mesh may coincide with the nodes and edges of the structured mesh. In these instances, we refer to the meshes as matching, as illustrated in Fig. 11. In such special cases, PDFEM coincides exactly with a classical FEM.

### 3.2. Computation of global properties with PDFEM on a periodic RVE

As previously mentioned, the Phantom domain Method is suitable for computing global properties of composite materials despite the lack of weak continuity at the interface. To substantiate this claim, we selected a periodic RVE 3D example featuring a fiber-like inclusion defined by a spline curve with a circular section, as depicted in Fig. 12. Notably, the left fiber is periodically repeated on the right.

The computation of global properties involves resolving six boundary value problems with periodic conditions applied to the faces of the Representative Volume Element (RVE) (Debordes, 1986), derived from a multi-scale modeling of composite materials with the definition of a periodic RVE, see Bensoussan et al. (1978), Sanchez-Palencia and Sanchez-Hubert (1992). We compared the results obtained from PDFEM calculations with a conformal Finite Element (FE) computation as a reference, using a mesh that was sufficiently fine for the purpose. The finite element computation proved intricate due to the need for multiple adjustments in the mesh construction; for example, the mesh must be adequately fine. We anticipate more complexities in constructing such meshes with an increased number of fibers or smaller thicknesses.

We opted for a second-order tetrahedron mesh with an element size of 0.05 (20 subdivisions on every edge of the RVE) and 30 subdivisions along each fiber, resulting in a 245,422 nodes mesh (see Fig. 13).

As shown in Fig 14, we observe an excellent linear convergence of the homogenized various elastic coefficients

To assess the capability of the method, we opted for a high contrast matrix/inclusions ratio of $c = 100$.

We present the computation time performed on a personal computer using a single core in Table 1. The Phantom domain Finite Element Method (PDFEM) proves its capacity to deliver swift solution times, even on consumer processors, compared to the size of the linear system. Specifically, without parallelization, less than 1 min have been required for an over one million degrees of freedom (dofs) problem using a Krylov-based solver (conjugate gradient with incomplete Cholesky preconditioner).

### 4. Conclusion

The Phantom Domain Finite Element Method (PDFEM) is an innovative approach for solving boundary value problems for heterogeneous material, such as thermal and elastic problems, inspired by Fictitious Domain Methods.

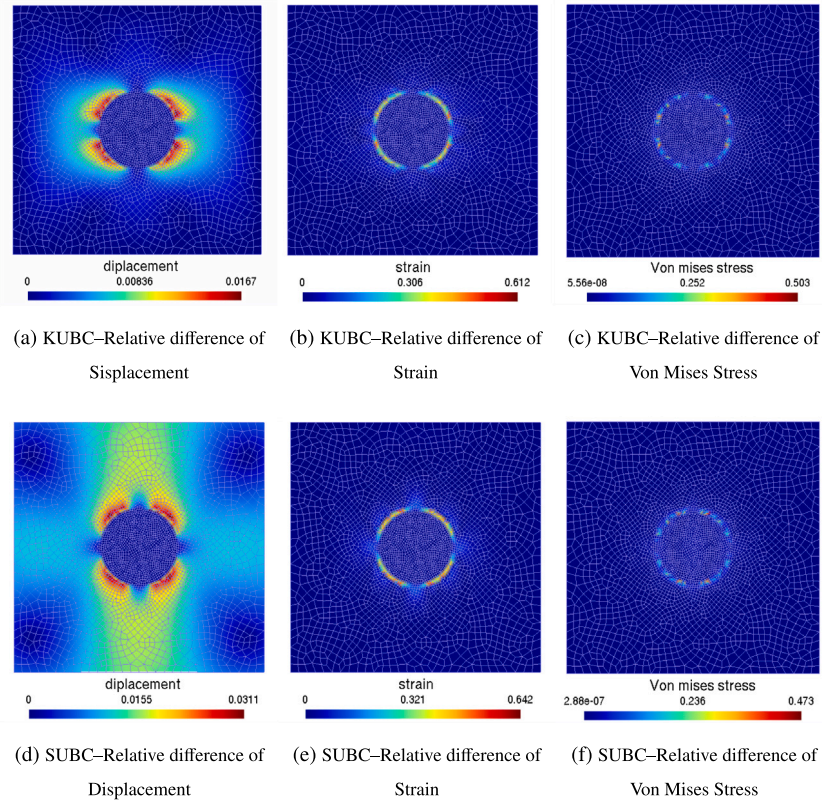

**Fig. 8.** PDFEM/FEM Relative difference for Displacement, Strain and Von Mises Stress for elastic boundary value problems KUBC (*a*), (*b*), (*c*) and SUBC (*d*), (*e*), (*f*). Case of one disk inclusion of diameter $d = 0.3$.

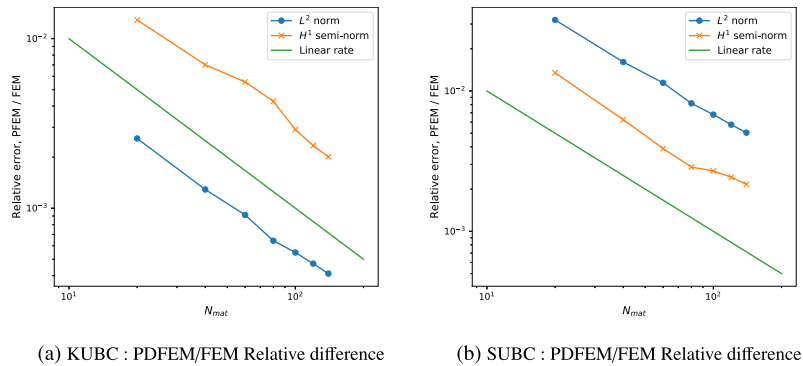

**Fig. 9.** PDFEM/FEM Relative difference of the displacement $L^2$ norm and $H^1$ semi-norm for a static elastic boundary value problems with respect to the mesh resolution. Case of one disk inclusion $d = 0.3$.

**Table 1**
Computation time for periodic homogenization using PDFEM (3D Elasticity).

| RVE parameter $n_r$ | 15 | 20 | 30 | 50 | 75 |
| --- | --- | --- | --- | --- | --- |
| Mesh Element size $h$ | 6.7 $10^{-2}$ | 5 $10^{-2}$ | 3.3 $10^{-2}$ | 2 $10^{-2}$ | 1.3 $10^{-2}$ |
| DOFs (in thousands) | 12K | 24K | 81K | 412K | 1 265K |
| Time (s) | 1.1 | 5 | 13 | 60.0 | 298 |

This method decomposes the energy functional of the variational problem into two distinct components: one representing the entire domain and the other corresponding to the inclusion. Unlike traditional Finite Element Methods (FEM), which require a conforming mesh, PDFEM employs independent meshes for these components.

By utilizing isoparametric elements on a structured mesh for the domain $\Omega$, we derive a substitution matrix that establishes a connection between the degrees of freedom (DOFs) of the inclusion mesh and those of the global structured mesh. This enables the implementation of FEM within the framework of a structured mesh for Representative Volume Elements (RVEs).

The substitution procedure can be interpreted as a smoothed pixelization of the inclusion within the structured mesh, requiring the inclusion mesh to be sufficiently refined relative to the structured mesh to ensure accurate approximations of the integral forms.

Numerical investigations performed on basic inclusions, such as disks and squares, demonstrate a linear convergence of relative errors when compared to reference FEM solutions. This convergence rate aligns with theoretical predictions based on linear interpolation. A comparison of PDFEM solutions with standard FEM solutions on the same mesh reveals that the error introduced by PDFEM is primarily concentrated around the matrix/inclusion interface due to the absence of

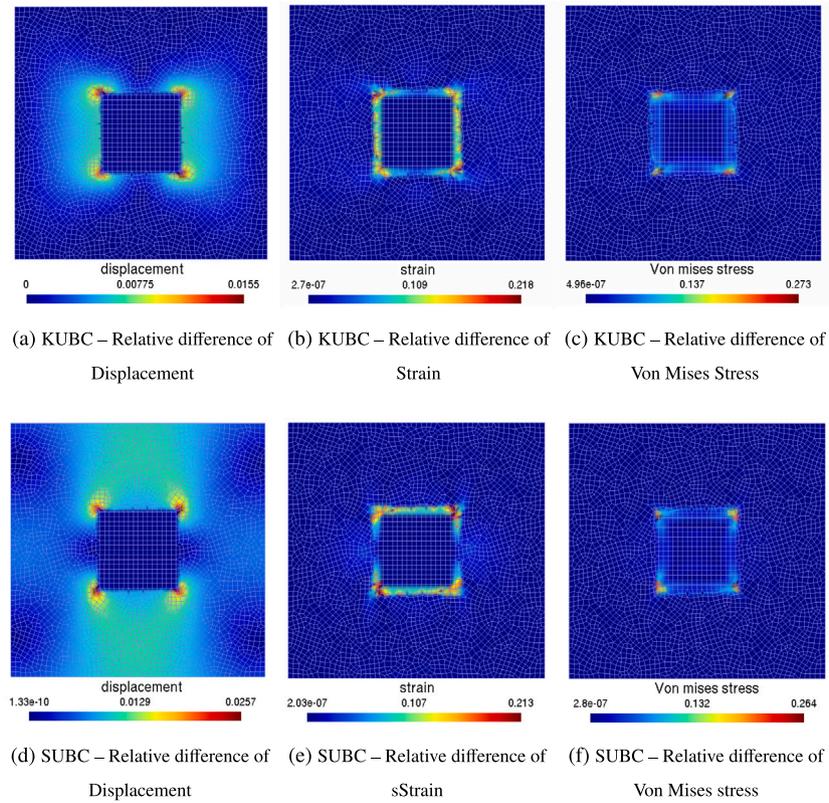

(a) KUBC – Relative difference of Displacement

(b) KUBC – Relative difference of Strain

(c) KUBC – Relative difference of Von Mises Stress

(d) SUBC – Relative difference of Displacement

(e) SUBC – Relative difference of sStrain

(f) SUBC – Relative difference of Von Mises stress

**Fig. 10.** PDFEM/FEM Relative difference for Displacement, Strain and Von Mises Stress for elastic boundary value problems KUBC (*a*), (*b*), (*c*) and SUBC (*d*), (*e*), (*f*). Case of one square inclusion of side length $l = 0.3$.

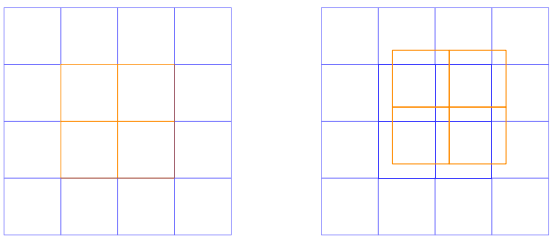

**Fig. 11.** A square inclusion in a structured mesh $\Omega$ of $N_{mat} = 4$ with matching meshes (left) and non-matching meshes (right).

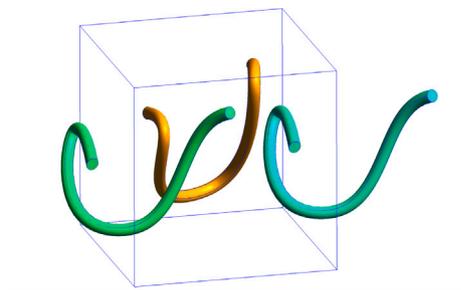

**Fig. 12.** A 3D periodic RVE with two fiber inclusions.

weak continuity, highlighting a current limitation of the method while still maintaining its effectiveness in computing the global properties of composite materials.

Through a simple example of a two-fiber inclusion periodic RVE, we demonstrate favorable convergence properties and computational efficiency. The primary advantage of PDFEM lies in its complete independence from mesh conformity, offering significant flexibility in handling intricate geometries, including curved fiber inclusions, which pose challenges for standard FEM methods. Moreover, PDFEM proves to be more time-efficient than conventional FEM for complex geometries, as it relies on a structured mesh, whereas FEM requires larger, often ill-conditioned matrices when dealing with intricate geometries.

The Phantom Domain Finite Element Method offers significant advantages in predicting the behavior of composite materials by computing stiffness and conductivity tensors for both mechanical and thermal problems. While further refinements are possible, the method already yields promising results, particularly in the assessment of global properties for complex geometries, such as randomly distributed fiber-reinforced composites.

Key advantages of our method include:

- Minimal computational complexity: The method avoids the need for enriched shape functions or complex interface treatments, making it easier to implement in standard finite element frameworks.
- Straightforward homogenization calculations: The approach is well-suited for evaluating effective material properties without requiring advanced meshing strategies.
- Energy-based formulation: By leveraging energetic principles, our method provides a robust yet computationally light framework for composite material modeling.

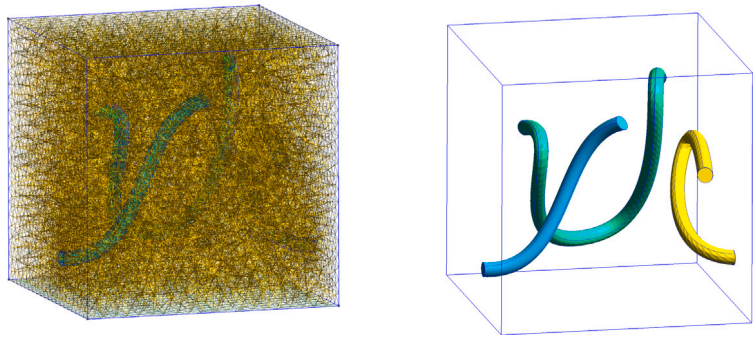

**Fig. 13.** Conformal periodic mesh with two fiber inclusions..

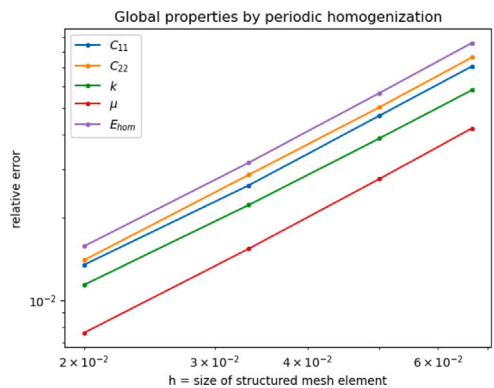

**Fig. 14.** Linear convergence PDFEM for periodic homogenization.

**CRediT authorship contribution statement**

**Tianlong He:** Writing – original draft, Software, Investigation. **Philippe Karamian-Surville:** Writing – review & editing, Validation, Supervision. **Daniel Choï:** Validation, Supervision, Software.

**Declaration of competing interest**

The authors declare the following financial interests/personal relationships which may be considered as potential competing interests: Philippe KARAMIAN-SURVILLE reports administrative support was provided by University of Caen Normandy. If there are other authors, they declare that they have no known competing financial interests or personal relationships that could have appeared to influence the work reported in this paper.

**Data availability**

No data was used for the research described in the article.